\documentclass{llncs}
\usepackage{graphicx}
\usepackage{amsmath}
\usepackage{amssymb}
\usepackage{enumerate}
\usepackage{varioref}
\usepackage{charter,eulervm}

\newcommand{\es}{\varnothing}

\title{\sc {The Black-and-White Coloring Problem on 
Distance-Hereditary Graphs 
and Strongly Chordal graphs}}

\author{
 Ton~Kloks\thanks{National 
Science Council of Taiwan Support Grant 
NSC~99--2218--E--007--016.}%
\inst{1}
\and 
 Sheung-Hung~Poon\inst{1}
\and 
 Feng-Ren~Tsai\inst{2} 
\and 
 Yue-Li~Wang\inst{3}
}
\institute{
 Department of Computer Science\\
 National Tsing Hua University,
 No.~101, Sec.~2, Kuang Fu Rd., Hsinchu, Taiwan\\
 {\tt spoon@cs.nthu.edu.tw}  
\and 
 Institute of Information Systems and Applications\\
 National Tsing Hua University,
 No.~101, Sec.~2, Kuang Fu Rd., Hsinchu, Taiwan\\
 {\tt mevernom@gmail.com}  
\and 
 Department of Information Management\\
 National Taiwan University of Science and Technology\\
 No.~43, Sec.~4, Keelung Rd., Taipei, 106, Taiwan\\
 {\tt ylwang@cs.ntust.edu.tw}
}

\begin{document}

\maketitle

\begin{abstract}
Given a graph $G$ and integers $b$ and $w$. 
The black-and-white coloring problem asks if there 
exist disjoint sets of vertices $B$ and $W$ with $|B|=b$ and $|W|=w$ 
such that no vertex in $B$ is adjacent to any vertex 
in $W$. In this paper we show that the problem is polynomial 
when restricted to cographs, distance-hereditary graphs, 
interval graphs and strongly chordal graphs. We show that the 
problem is NP-complete on splitgraphs. 
\end{abstract}

\section{Introduction}

\begin{definition}
Let $G=(V,E)$ be a graph and let $b$ and $w$ be two integers. 
A black-and-white coloring of $G$ colors $b$ vertices black 
and $w$ vertices white such that no black vertex is adjacent to 
any white vertex. 
\end{definition}

In other words, the black-and-white coloring problem 
asks for a complete bipartite subgraph $M$ in the 
complement $\Bar{G}$ of $G$ with $b$ and $w$ vertices in the 
two color classes of $M$. 

The black-and-white coloring problem is NP-complete 
for graphs in general~\cite{kn:hansen}. 
That paper also shows that 
the problem can be solved for trees in $O(n^3)$ time. 
In a recent paper~\cite{kn:berend} the worst-case timebound 
for an algorithm on  
trees was improved to $O(n^2 \log^3 n)$ time~\cite{kn:berend}. 
The paper~\cite{kn:berend} mentions, among other things, 
a manuscript by Kobler, 
{\em et al.\/}, which shows that the problem can be solved in 
polynomial time for graphs of bounded treewidth. 

\medskip 

In this paper we investigate the complexity of the problem 
for some graph classes. We start our analysis for  
the class of cographs.  

A $P_4$ is a path with four vertices. 
 
\begin{definition}[\cite{kn:corneil}]
A graph is a cograph if it has no induced $P_4$. 
\end{definition}

There are various characterizations of cographs. 
For algorithmic purposes the following characterization 
is suitable. 

\begin{theorem}
A graph is a cographs if and only if every induced 
subgraph $H$ is disconnected or the complement $\Bar{H}$ 
is disconnected. 
\end{theorem}

It follows that a cograph has a tree decomposition which 
is called a cotree. A cotree  
is a pair $(T,f)$ comprising  
a rooted binary tree $T$ together with  
a bijection $f$ from the vertices 
of the graph to the leaves of the tree. Each internal node 
of $T$, including the root, has a label $\otimes$ or $\oplus$. 
The $\otimes$ operation is called a 
join operation, and it makes every vertex that is mapped to 
a leaf in the left subtree adjacent to every vertex that is 
mapped to a leaf in the right subtree. The operator  
$\oplus$ is called a union operation. In that 
case the graph is the union of the graphs defined by the 
left - and right subtree. A cotree decomposition can be 
obtained in linear time~\cite{kn:corneil2}.  
  
\section{Black-and-white colorings of cographs}
\label{section cographs}

In this section we show that the black-and-white coloring 
problem can be solved in polynomial time for cographs. 

\begin{lemma}
\label{cg}
There exists an $O(n^5)$ algorithm which solves the 
black-and-white coloring problem on cographs. 
\end{lemma} 
\begin{proof}
Let $\gamma_G(b,w)$ be a boolean variable which indicates if the 
graph $G$ has a black-and-white coloring with $b$ black and 
$w$ white vertices. 
Obviously, we have 
\begin{equation}
\label{eqn}
\gamma_G(b,0)= 
\begin{cases} 
\mathrm{true} & \quad\text{if $0 \leq b \leq n$,}\\
\mathrm{false} & \quad\text{otherwise,}
\end{cases}
\end{equation}
where $n$ is the number of vertices in $G$. A similar 
formula holds for $\gamma_G(0,w)$. 
The algorithm uses dynamic programming on the cotree and 
it derives $\gamma_G(b,w)$ for every node from a 
table of values stored at the children as 
follows. 

\smallskip 

When $G$ has only one vertex then we have 
\[\gamma_G(b,w)=
\begin{cases} 
\mathrm{true} & \quad\text{if $(b,w) \in 
\{\; (0,0), \;(1,0), \;(0,1)\;\}$,}\\ 
\mathrm{false} & \quad\text{in all other cases.} 
\end{cases}\]

\smallskip 

Assume that $G$ is the join of two cographs $G_1$ and $G_2$. 
Let $n_i$ be the number of vertices of $G_i$. Then $n=n_1+n_2$.  
If there is a black-and-white coloring for $G$ with at least one 
black vertex and at least one white vertex then all the black 
and white vertices must be contained in the same graph $G_i$. 
It follows that, when $b\geq 1$ and $w \geq 1$,   
\begin{eqnarray*}
\gamma_G(b,w)&=&\text{true} \quad\text{if and only if}\\ 
&& \gamma_{G_1}(b,w)=\text{true} \quad\text{or}\quad  
\gamma_{G_2}(b,w)=\text{true}.
\end{eqnarray*}
The cases where $b=0$ or $w=0$ follow from Equation~\ref{eqn}. 
 
\smallskip 

Finally assume that $G$ is the union of two cographs $G_1$ and 
$G_2$. Let $n_i$ be the number of vertices in $G_i$ and let 
$n=n_1+n_2$ be the number of vertices in $G$. 
Then we have 
\begin{eqnarray*}
\gamma_G(b,w) &= &\text{true} \quad\text{if and only if}\\ 
&& \exists_{k}\;\exists_{\ell} \;\; 
\gamma_{G_1}(k,\ell)=\text{true} \quad\text{and}\quad 
\gamma_{G_2}(b-k,w-\ell)=\text{true}.  
\end{eqnarray*}

\smallskip 

A table containing the boolean values $\gamma_G(b,w)$ 
has $n^2$ entries. By the formulas 
above, each entry can be computed in $O(n^2)$ time. Thus a 
complete table for each node in the cotree can be obtained in 
$O(n^4)$ time. Since a cotree has $O(n)$ nodes, this algorithm 
can be implemented to run in $O(n^5)$ time. 
\qed\end{proof}

The following theorem improves the timebound. 

\begin{theorem}
There exists an $O(n^3)$ algorithm which solves the 
black-and-white coloring problem on cographs.
\end{theorem}
\begin{proof}
Let $f_G(b)$ be the maximum number of white 
vertices in a black-and-white coloring of $G$ 
with $b$ black vertices. 
We prove that the function $f_G$ can be computed in 
$O(n^3)$ time for cographs. 

Let $G$ be a cograph with $n$ vertices. 
We write $f$ instead of $f_G$. 
By convention, 
\[f(b)=0 \quad\text{when $b < 0$ or $b > n$.}\] 
  
\smallskip 

Assume that $G$ has one vertex. Then 
\[f(b) = 
\begin{cases}
1 & \text{if $b=0$} \\
0 & \text{in all other cases.}
\end{cases}\] 

\smallskip 

Assume that $G$ is the join of two cographs $G_1$ 
and $G_2$. We write $f_i$ instead of $f_{G_i}$, for 
$i \in \{1,2\}$. 
We have that $f(0)=n$, where $n$ is the number 
of vertices in $G$. 
When $b > 0$ we have   
\[f(b) = \max \; \{\;f_1(b),\;f_2(b)\;\}.\] 
 
\smallskip 

Assume that $G$ is the union of two cographs $G_1$ and $G_2$. 
Then 
\[f(b)= \max_{0 \leq k \leq b} \; f_1(k)+f_2(b-k).\]

\smallskip 

A cotree $T$ has $O(n)$ nodes and it can be computed 
in linear time~\cite{kn:corneil}. Consider a node $i$ in $T$. Let 
$G_i$ be the subgraph of $G$ induced by the vertices that are 
mapped to leaves in the subtree rooted at $i$. 
By the previous observations, the function $f_i$ for the 
graph $G_i$ can be computed in $O(n^2)$ time. 
Since $T$ has $O(n)$ nodes this proves the theorem. 
\qed\end{proof}

\subsection{Threshold graphs}
\label{section threshold}

A subclass of the class of cographs are the threshold graphs. 

\begin{definition}[\cite{kn:chvatal}]
A graph $G=(V,E)$ is a threshold graph if there is a real number 
$T$ and a real number $w(x)$ for every vertex $x \in V$ such that 
a subset $S \subseteq V$ is an independent set if and only if 
\[\sum_{x \in S} w(x) \geq T.\]
\end{definition}

There are many ways to characterize threshold graphs~\cite{kn:mahadev}. 
For example, 
a graph is a threshold graph if it has no induced $P_4$, $C_4$ or 
$2K_2$. 
\begin{figure}
\setlength{\unitlength}{1.35pt}
\begin{center}
\begin{picture}(140,20)
\thicklines
\put(0,0){\circle*{3.0}}
\put(20,0){\circle*{3.0}}
\put(0,20){\circle*{3.0}}
\put(20,20){\circle*{3.0}}
\put(0,0){\line(1,0){20}}
\put(0,0){\line(0,1){20}}
\put(20,20){\line(-1,0){20}}
\put(20,20){\line(0,-1){20}}

\put(60,0){\circle*{3.0}}
\put(80,0){\circle*{3.0}}
\put(60,20){\circle*{3.0}}
\put(80,20){\circle*{3.0}}
\put(60,0){\line(0,1){20}}
\put(80,20){\line(-1,0){20}}
\put(80,20){\line(0,-1){20}}

\put(120,0){\circle*{3.0}}
\put(140,0){\circle*{3.0}}
\put(120,20){\circle*{3.0}}
\put(140,20){\circle*{3.0}}
\put(120,0){\line(1,0){20}}
\put(120,20){\line(1,0){20}}
\end{picture}
\caption{A graph is a threshold graph if it has 
no induced $C_4$, $P_4$ or $2K_2$.}
\label{fig threshold}
\end{center}
\end{figure}
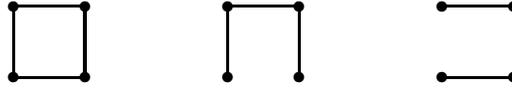

Another characterization is that a graph is a threshold graph 
if every induced subgraph has a universal vertex or an isolated 
vertex~\cite[Theorem 1]{kn:chvatal}.  
In~\cite[Corollary 1B]{kn:chvatal} appears also the 
following characterization. A graph $G=(V,E)$ 
is a threshold graph if and only if there is a partition 
of $V$ into two sets $A$ and $B$, of which one is possibly empty, 
such that 
\begin{enumerate}[\rm 1.]
\item $A$ induces a clique, 
\item $B$ induces an independent set, and 
\item there is an ordering $b_1,\ldots,b_k$ of the 
vertices in $B$ such that 
\[N(b_1) \subseteq \ldots \subseteq N(b_k).\]
\end{enumerate}

\medskip 

We use the notation $N[x]$ to denote the closed neighborhood 
of a vertex $x$. Thus $N[x]=N(x) \cup \{x\}$. 

\begin{theorem}
\label{thm threshold}
There exists a linear-time algorithm which, 
given a threshold graph $G$ and integers $b$ and $w$, 
decides if there is a black-and-white coloring of $G$ 
with $b$ vertices colored black and $w$ vertices colored white. 
\end{theorem}
\begin{proof}
Let $x_1,\ldots,x_n$ be an ordering of the vertices in $G$ 
such that for all $i < n$ 
\begin{enumerate}[\rm (a)]
\item $N(x_i) \subseteq N(x_{i+1})$ if $x_i$ and $x_{i+1}$ are 
not adjacent, and 
\item $N[x_i] \subseteq N[x_{i+1}]$ if $x_i$ and $x_{i+1}$ are adjacent. 
\end{enumerate}

Assume that there exists a black-and-white coloring 
which colors $b$ vertices black and $w$ vertices 
white. 
Assume that there is an index $k \leq b+w$ such that 
$x_k$ is uncolored. Then there exists an 
index $\ell > b+w$ such that $x_{\ell}$ is colored black or white. 
Then we may color $x_k$ with the color of $x_{\ell}$ and 
uncolor $x_{\ell}$ instead. Thus we may assume 
that there exists a coloring 
such that $x_1,\ldots,x_{b+w}$ are colored and  
all other vertices are uncolored. 

\smallskip 

Assume that $b \leq w$. 
We prove that there exists a coloring $f$ such that 
\[f(x_i)= 
\begin{cases}
\text{black} & \text{if $1 \leq i \leq b$, and}\\
\text{white} & \text{if $b+1 \leq i \leq b+w$.} 
\end{cases}\]

We may assume that $b \geq 1$ and that $w \geq 1$. 
Assume that $x_i$ is adjacent to $x_j$ for some 
$i \leq b < j$. Then 
\[\{x_j,\ldots,x_{b+w}\}  \subseteq  N(x_i)
\quad\text{and}\quad \{x_i,\ldots,x_j\}  \subseteq  N[x_j].\]
Thus all vertices in 
\[\{x_i,\ldots,x_{b+w}\}\]  
are the same color. If they are all black then 
there are at least $w+1$ black vertices in the coloring, 
which contradicts $b \leq w$. If they are all white  
then we have at least $w+1$ white vertices, which  
is a contradiction as well. 
Thus no two vertices $x_i$ and 
$x_j$ with $i \leq b < j$ are adjacent, which proves 
that the coloring above is valid. 

This proves the theorem, since an algorithm 
only needs to check if $x_b$ is adjacent to $x_{b+w}$ or not. 
\qed\end{proof}

\subsection{Difference graphs}

\begin{definition}[\cite{kn:hammer}]
A graph $G=(V,E)$ is a difference graph if there exists a 
positive real number 
$T$ and a real number $w(x)$ for every vertex $x \in V$ 
such that  
$w(x) \leq T$ for every $x \in V$ and such that 
for any pair of vertices $x$ and $y$ 
\[\{x,y\} \in E \quad\text{if and only if}\quad 
|w(x)-w(y)| \geq T.\]
\end{definition}

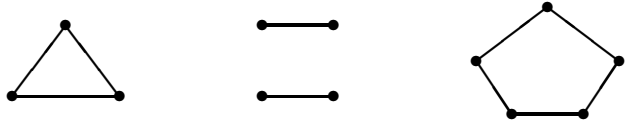
\begin{figure}
\setlength{\unitlength}{1.35pt}
\begin{center}
\begin{picture}(175,30)
\thicklines
\put(0,5){\circle*{3.0}}
\put(30,5){\circle*{3.0}}
\put(15,25){\circle*{3.0}}
\put(0,5){\line(1,0){30}}
\put(0,5){\line(3,4){15}}
\put(30,5){\line(-3,4){15}}

\put(70,5){\circle*{3.0}}
\put(90,5){\circle*{3.0}}
\put(70,25){\circle*{3.0}}
\put(90,25){\circle*{3.0}}
\put(70,5){\line(1,0){20}}
\put(70,25){\line(1,0){20}}

\put(130,15){\circle*{3.0}}
\put(140,0){\circle*{3.0}}
\put(150,30){\circle*{3.0}}
\put(160,0){\circle*{3.0}}
\put(170,15){\circle*{3.0}}
\put(130,15){\line(2,-3){10}}
\put(130,15){\line(4,3){20}}
\put(140,0){\line(1,0){20}}
\put(160,0){\line(2,3){10}}
\put(170,15){\line(-4,3){20}}
\end{picture}
\caption{A graph is a difference graph if it has
no induced triangle, $2K_2$ or $C_5$.}
\label{fig difference}
\end{center}
\end{figure}

Difference graphs are sometimes called 
chain graphs~\cite{kn:yannakakis}. 

\medskip 

Difference graphs can be characterized in 
many ways~\cite{kn:hammer}. For example, a graph 
is a difference graph if and only if it has 
no induced 
$K_3$, $2K_2$ or $C_5$~\cite[Proposition~2.6]{kn:hammer}. 
Difference graphs are bipartite. Let $X$ and $Y$ 
be a partition of $V$ into two color classes. 
Then the graph obtained 
by making a clique of $X$ is a threshold graph 
and this property characterizes difference 
graphs~\cite[Lemma~2.1]{kn:hammer}. 

\begin{theorem}
\label{thm difference}
There exists a linear-time algorithm which, given a 
difference graph $G$ and integers $b$ and $w$, decides if there 
is a black-and-white coloring of $G$ with $b$ black vertices and 
$w$ white vertices. 
\end{theorem}
\begin{proof}
An argument, similar to the one given in~Theorem~\ref{thm threshold}, 
provides the proof. 
\qed\end{proof}
    
\section{Distance-hereditary graphs}
\label{section DH}

\begin{definition}[\cite{kn:howorka}]
A graph $G$ is distance hereditary if for every pair of 
nonadjacent vertices 
$x$ and $y$ and for every connected induced subgraph $H$ of $G$ 
which contains $x$ and $y$, the distance between $x$ and $y$ 
in $H$ is the same as the distance between $x$ and $y$ in $G$. 
\end{definition}

In other words, a graph $G$ is distance hereditary if 
for every nonadjacent pair 
$x$ and $y$ of vertices, 
all chordless paths between $x$ and $y$ 
in $G$ have the same length. 

\medskip 

There are various characterizations of distance-hereditary graphs. 
One of them states that a graph is distance hereditary if and 
only if it has no induced house, hole, domino or 
gem~\cite{kn:bandelt,kn:howorka}. 
Distance-hereditary graphs are also characterized 
by 
the property that every induced subgraph has either an isolated 
vertex, or a pendant vertex, or a true or false 
twin~\cite{kn:bandelt}. 

\begin{figure}
\setlength{\unitlength}{1.35pt}
\begin{center}
\begin{picture}(230,30)
\thicklines
\put(0,0){\circle*{3.0}} 
\put(20,0){\circle*{3.0}}
\put(0,20){\circle*{3.0}} 
\put(20,20){\circle*{3.0}}
\put(10,30){\circle*{3.0}} 
\put(0,0){\line(1,0){20}}
\put(0,0){\line(0,1){20}} 
\put(20,20){\line(-1,0){20}}
\put(20,20){\line(0,-1){20}} 
\put(10,30){\line(-1,-1){10}}
\put(10,30){\line(1,-1){10}}

\put(50,15){\circle*{3.0}} 
\put(65,0){\circle*{3.0}}
\put(65,30){\circle*{3.0}} 
\put(85,0){\circle*{3.0}}
\put(100,15){\circle*{3.0}} 
\put(50,15){\line(1,-1){15}}
\put(50,15){\line(1,1){15}} 
\put(65,0){\line(1,0){20}}
\put(85,0){\line(1,1){15}} 
\qbezier[25](65,30)(90,45)(100,15)

\put(130,0){\circle*{3.0}}
\put(150,0){\circle*{3.0}}
\put(130,15){\circle*{3.0}}
\put(150,15){\circle*{3.0}}
\put(130,30){\circle*{3.0}}
\put(150,30){\circle*{3.0}}
\put(130,0){\line(0,1){30}} 
\put(150,0){\line(0,1){30}}
\put(130,0){\line(1,0){20}}
\put(130,15){\line(1,0){20}}
\put(130,30){\line(1,0){20}}

\put(180,15){\circle*{3.0}} 
\put(195,30){\circle*{3.0}}
\put(205,0){\circle*{3.0}} 
\put(215,30){\circle*{3.0}}
\put(230,15){\circle*{3.0}} 
\put(205,0){\line(-5,3){25}}
\put(205,0){\line(-1,3){10}}
\put(205,0){\line(1,3){10}}
\put(205,0){\line(5,3){25}} 
\put(180,15){\line(1,1){15}}
\put(195,30){\line(1,0){20}} 
\put(215,30){\line(1,-1){15}}
\end{picture}
\caption{A graph is distance hereditary if it has no induced
house, hole, domino or gem.}
\label{HHDG}
\end{center}
\end{figure}
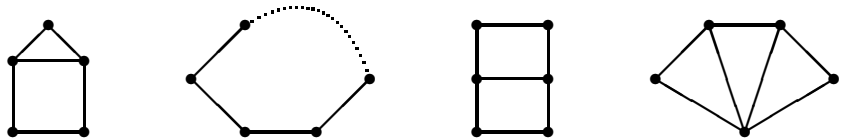

\medskip 

Distance-hereditary graphs are the graphs of 
rankwidth one. This implies that they have a 
special decomposition tree which we describe next. 
 
\medskip 

A decomposition tree for a graph $G=(V,E)$ is a pair $(T,f)$ 
consisting of   
a rooted binary tree $T$ and a bijection $f$ from  
$V$ to the leaves of $T$. 

\medskip 

When $G$ is distance hereditary it has a 
decomposition 
tree $(T,f)$ with the following 
three properties~\cite{kn:chang}. 

Consider an edge $e=\{p,c\}$ in $T$ where $p$ is the parent of $c$. 
Let $W_e \subset V$ be the set of vertices of $G$ 
that are mapped by $f$ to the leaves 
in the subtree rooted at $c$. Let $Q_e \subseteq W_e$ be the 
set of vertices in $W_e$ that have neighbors in $G-W_e$. 
The set $Q_e$ is called the twinset of $e$. 
The first property is that the subgraph of 
$G$ induced by $Q_e$ is a cograph for every 
edge $e$ in $T$. 

\medskip 

Consider an internal vertex $p$ in $T$. 
Let $c_1$ and 
$c_2$ be the two children of $p$. Let $e_1=\{p,c_1\}$ and let 
$e_2=\{p,c_2\}$. Let $Q_1$ and $Q_2$ be the twinsets of $e_1$ 
and $e_2$. The second property is that 
there is a join- or a union-operation 
between $Q_1$ and $Q_2$. Thus every vertex of $Q_1$ 
has the same neighbors in $Q_2$. 

\medskip 

Let $p$ be an internal vertex of $T$ which is not the root. 
Let $e$ be the line that connects $p$ with its parent. 
Let $Q_e$ be the twinset of $e$.  
Let $c_1$ and $c_2$ be the two children of $p$ in $T$.  
Let $e_1=\{p,c_1\}$ and let $e_2=\{p,c_2\}$. 
Let $Q_i$ be the twinset of $e_i$, for $i \in \{1,2\}$. 
The third, and final, property is that 
\[Q_e=Q_1 \quad\text{or}\quad Q_e=Q_2 \quad\text{or}\quad 
Q_e=Q_1 \cup Q_2.\] 

\medskip 

When $G$ is distance hereditary then a tree-decomposition for 
$G$ with the three properties described above 
can be obtained in linear time~\cite{kn:chang}. 

\medskip 

Notice that the first property is a consequence of the other two. 
As an example, notice that cographs are distance hereditary. A cotree 
is a decomposition tree for a cograph with the three properties 
mentioned above.  
  
\begin{theorem}
There exists a polynomial-time algorithm that solves the 
black-and-white coloring problem on distance-hereditary graphs. 
\end{theorem}
\begin{proof}
Let $(T,f)$ be a tree-decomposition for $G$ 
which satisfies the properties mentioned above. 

\smallskip 

Define a boolean variable 
\[\gamma_e(b,w,b^{\prime},w^{\prime})\]  
for a subgraph induced by a branch rooted at a line $e$ of $T$. 
This variable is true if there exists a black-and-white 
coloring of the subgraph, induced by the vertices 
that are mapped to the leaves in the branch,  
with $b$ black vertices and $w$ white vertices,  
such that $b^{\prime}$ black vertices and 
$w^{\prime}$ white vertices are contained in the 
twinset of the branch. 

\smallskip 

First assume that $e=\{p,c\}$ is an edge of $T$ which connects a leaf 
$c$ with its parent. Let $Q$ be the twinset of $e$, that is, 
$Q=\es$ or $Q=\{c\}$. 
Then we have 

\begin{enumerate}[\rm (a)]
\item $\gamma_e(0,0,0,0)=\text{true}$, 
\item $\gamma_e(1,0,0,0)= \gamma_e(0,1,0,0)=\text{true}$ if $Q=\es$, 
\item $\gamma_e(1,0,1,0)=\gamma_e(0,1,0,1)=\text{true}$ if $Q=\{c\}$, and 
\item $\gamma_e(b,w,b^{\prime},w^{\prime})=\text{false}$ in all 
other cases.
\end{enumerate}

\smallskip 

Consider a node $p$ in $T$ with two children $c_1$ and $c_2$.  
Let $e_i=\{p,c_i\}$ for $i \in \{1,2\}$. Let $Q_1$ and 
$Q_2$ be the 
twinsets of $e_1$ and $e_2$. 
If $p$ is not the root then let 
$Q$ be the twinset for the line that connects $p$ with its 
parent. 
We consider the following cases. 
First assume that there is a join between $Q_1$ and $Q_2$ and 
that $Q=Q_1 \cup Q_2$. 
Consider the case where there are no white vertices 
in the twinset $Q$. 
Then we have, for all $p$, $b$, $w$ 
\[\gamma_e(b,w,p,0)=\text{true}\] 
if and only if there exists partitions $p=p_1+p_2$, $w=w_1+w_2$ and 
$b=b_1+b_2$ such that 
\[\gamma_{e_1}(b_1,w_1,p_1,0)=\text{true} \quad\text{and}\quad  
\gamma_{e_2}(b_2,w_2,p_2,0)=\text{true}.\] 
A similar formula holds for the case where there are no black 
vertices in the twinset. 

Next, consider black-and-white colorings 
where there is at least one black, and at least one 
white vertex in the twinset $Q$. 
Then we have that all the black and white vertices of $Q$ must be in 
one of $Q_1$ and $Q_2$. 
In that case we have, for all $b$, $w$,  
and for all $p>0$ and $q>0$ 
\[\gamma_e(b,w,p,q) = \text{true}\] 
if and only if there exist partitions 
$b=b_1+b_2$, $w=w_1+w_2$ such that  
\begin{align*}
(\gamma_{e_1}(b_1,w_1,p,q)=\text{true} &\quad\text{and} \quad  
\gamma_{e_2}(b_2,w_2,0,0)=\text{true}) \quad\text{or} \\
(\gamma_{e_1}(b_1,w_1,0,0)=\text{true} &\quad\text{and}\quad   
\gamma_{e_2}(b_2,w_2,p,q)=\text{true}).   
\end{align*}

\smallskip 

In the second case we assume that there is a join between 
$Q_1$ and $Q_2$ and that $Q=Q_1$. 
Obviously, we obtain the same formulas as above, except that 
the numbers of black and white vertices in the twinset $Q$ 
are copied from those numbers in $Q_1$. 
Thus we obtain that  
\[\gamma_e(b,w,p,q)=\text{true}\] 
if and only if there exist partitions $b=b_1+b_2$, 
and $w=w_1+w_2$ such that the following hold.  
\begin{eqnarray*}
\lefteqn{\text{If $p>0$ and $q>0$:}}\\ 
&& \gamma_{e_1}(b_1,w_1,p,q)=\text{true} \quad\text{and}\quad 
 \gamma_{e_2}(b_2,w_2,0,0) =\text{true} \\
\lefteqn{\text{if $p>0$ and $q=0$:}}\\
&& \gamma_{e_1}(b_1,w_1,p,0)=\text{true} \quad\text{and}\quad   
\exists_{p_2}\; \gamma_{e_2}(b_2,w_2,p_2,0)=\text{true}\\
\lefteqn{\text{if $p=0$ and $q>0$:}}\\
&&\gamma_{e_1}(b_1,w_1,0,q)=\text{true} \quad\text{and}\quad     
\exists_{q_2}\; \gamma_{e_2}(b_2,w_2,0,q_2)=\text{true}\\
\lefteqn{\text{and, if $p=0$ and $q=0$:}}\\ 
&&\gamma_{e_1}(b_1,w_1,0,0)=\text{true} \quad\text{and}\quad   
\exists_{p_2}\;\exists_{q_2}\; 
\gamma_{e_2}(b_2,w_2,p_2,q_2)=\text{true.}
\end{eqnarray*} 

\smallskip 

Now assume that there is a union between $Q_1$ and $Q_2$. 
First assume that $Q=Q_1 \cup Q_2$. 
Then we have for all $b$, $w$, $p$ and $q$,     
\[\gamma_e(b,w,p,q)=\text{true}\]
if and only if
there exist partitions $b=b_1+b_2$, 
$w=w_1+w_2$, $p=p_1+p_2$ and $q=q_1+q_2$ 
such that 
\[\gamma_{e_1}(b_1,w_1,p_1,q_1)=\text{true} \quad\text{and}\quad  
\gamma_{e_2}(b_2,w_2,p_2,q_2)=\text{true}.\] 

\smallskip 

Finally, assume that there is a union between $Q_1$ and 
$Q_2$ and that $Q=Q_1$. 
Then 
\[\gamma_e(b,w,p,q)=\text{true}\] 
if and only if there exist partitions 
$b=b_1+b_2$ and $w=w_1+w_2$, such that  
\[\gamma_{e_1}(b_1,w_1,p,q)=\text{true} \quad\text{and}\quad  
\exists_{p_2}\;\exists_{q_2}\; 
\gamma_{e_2}(b_2,w_2,p_2,q_2)=\text{true}.\]

\smallskip 

By symmetry, the remaining cases are similar. 

\smallskip 

When $p$ is the root, then the twinset $Q$ is not defined. 
To get around this obstacle we may simply add an edge $\Hat{e}$ 
in the tree adjacent to $p$ and define the twinset 
$Q$ for this edge, arbitrarily, as $Q=Q_1 \cup Q_2$, 
or $Q=Q_1$, or $Q=Q_2$. There exists a black-and-white 
coloring of $G$ with $b$ black and $w$ white vertices if there 
are $p$ and $q$ such that 
\[\gamma_{\Hat{e}}(b,w,p,q)=\text{true}.\]

A table consists of $O(n^4)$ entries for values of 
$b$, $w$, $p$ and $q$ ranging from $0$ up to $n$. 
For each node in the tree-decomposition, 
the value of each entry in the table can be 
computed in $O(n^8)$ time from the tables that are stored 
at the two children of the node. Therefore, a table at each node 
can be computed in $O(n^{12})$ time. Since the 
tree-decomposition has $O(n)$ nodes, this gives an 
upperbound of $O(n^{13})$ for solving the black-and-white 
coloring problem on distance-hereditary graphs.  
\qed\end{proof}

\section{Interval graphs}
\label{section interval}

In this section we show that there is an efficient algorithm 
to solve the black-and-white coloring problem on interval graphs. 

\begin{definition}[\cite{kn:lekkerkerker}]
A graph $G$ is an interval graph if it is the intersection 
graph of a collection of intervals on the real line. 
\end{definition}

There are various characterizations of interval graphs. 
For example, a graph is an interval graph if and only if 
it is chordal and it has no asteroidal triple. 
Also, a graph is an interval 
graph if and only if it has no $C_4$ and the complement 
$\Bar{G}$ has a transitive orientation~\cite{kn:gilmore}. 

\medskip 

For our purposes the following characterization of 
interval graphs is suitable. 

\begin{theorem}[\cite{kn:gilmore}]
A graph $G$ is an interval graph if and only if there 
is a linear ordering $L$ of its maximal cliques such that 
for every vertex, the maximal cliques that contain that 
vertex are consecutive in $L$. 
\end{theorem}

Interval graphs can be recognized in linear time. When 
$G$ is an interval graph then $G$ is chordal and so it has 
at most $n$ maximal cliques. A linear ordering of the maximal 
cliques can be obtained in $O(n^2)$ time~\cite{kn:booth}. 

\begin{theorem}
There exists an $O(n^6)$ algorithm that solves 
the black-and-white coloring problem on interval graphs. 
\end{theorem}
\begin{proof}
Let $[C_1,\ldots,C_t]$ be a linear ordering of the maximal 
cliques of an interval graph $G=(V,E)$ such that for every 
vertex $x$,  
the maximal cliques that contain $x$ appear consecutively 
in this ordering. 

\smallskip 

Consider a black-and-white coloring of $G$. First assume 
that the first clique $C_1$ contains no black or white 
vertices. Then we may remove the vertices that appear in $C_1$ 
from the graph and consider a black-and-white 
coloring of the vertices in cliques of the linear ordering 
\[[C_2^{\ast},\ldots,C_t^{\ast}],  
\quad\text{where, for $i > 1$,}\quad C_i^{\ast}=C_i \setminus C_1.\]
  
Now assume that $C_1$ contains some black vertices. 
Then, obviously, $C_1$ contains no white vertices. 
Let $i$ be the maximal index such that all the cliques 
$C_{\ell}$ with $1 \leq \ell \leq i$ contain no white vertices. 
Remove 
all the vertices that appear in $C_1,\ldots,C_i$ from the 
remaining cliques and consider the ordering 
\[ [C_{i+1}^{\ast},\ldots,C_t^{\ast}] \quad\text{where, 
for $\ell > i$,}\quad   
C_{\ell}^{\ast}=C_{\ell} \setminus \bigcup_{k=1}^i C_k.\] 
Then we may take an arbitrary black-and-white coloring 
of the graph induced by the vertices $\cup_{\ell=i+1}^t C_{\ell}^{\ast}$ 
and color an arbitrary number of vertices in 
$\cup_{\ell=1}^i C_{\ell}$
black. 

\smallskip 

For this purpose define, for $p \leq  q$, 
\[X_{p,q}=\{\; x \in V \;|\; x \in C_k \quad\text{if and only if}\quad 
p \leq k \leq q \;\}.\] 
Thus $X_{p,q}$ consists of the vertices of which the 
indices of the first and the last clique that contain the vertex 
are both in the interval $[p,q]$. 

\smallskip 

For $i \geq 1$ let $G_i$ be the graph with vertices in 
\[\bigcup_{k=i}^t C_k^{i}, \quad\text{where, for $k \geq i$,}\quad  
C_k^{i}=C_k \setminus 
\bigcup_{\ell=1}^{i-1} C_{\ell}.\]

\smallskip 

The algorithm keeps a table with entries $b,w \in \{1,\ldots,n\}$ 
and the boolean value $\gamma_i(b,w)$ which is true if and only 
if there exists a black-and-white coloring of $G_i$ 
with $b$ black vertices and $w$ white vertices. 
Then we have, for $i=1,\ldots,t$,   
\begin{eqnarray*}
\gamma_i(b,w) &=& \text{true} \quad\text{if and only if}
\quad \exists_{j \geq i} \; \exists_k \;\; 0 \leq k \leq |X_{i,j}| 
\quad\text{and} \\   
&& 
\begin{cases} 
\text{$(b,w) \in \{(k,0),(0,k)\}$} 
& \quad \text{if $j =t$, and } \\
\gamma_{j+1}(b-k,w) \;\text{or}\; \gamma_{j+1}(b,w-k)  
& \quad \text{if $j < t$.}
\end{cases}
\end{eqnarray*}

\smallskip 

To implement this algorithm one needs to 
compute the cardinalities $|X_{p,q}|$. 
Initialize $|X_{p,q}|=0$. We assume that we have, for each 
vertex $x$, the index $F(x)$ of the first clique 
that contains $x$ and the index $L(x)$ 
of the last clique   
that contains $x$. Consider the vertices one by one. 
For a vertex $x$, add one to 
$|X_{p,q}|$ for all $p \leq F(x)$ and all $q \geq L(x)$. 
For each vertex $x$ we need to update $O(n^2)$ cardinalities 
$|X_{p,q}|$. Thus computing all cardinalities $|X_{p,q}|$ can be 
done in $O(n^3)$ time. 

\smallskip 

For each $i=1,\ldots,t$, the table for $G_i$ contains 
$O(n^2)$ boolean values $\gamma_i(b,w)$. For the computation 
of each $\gamma_i(b,w)$ the algorithm searches the tables 
of $G_j$ for all $j > i$. Thus the computation of $\gamma_i(b,w)$ 
takes $O(n^3)$ time. Thus the full table for $G_i$ can be 
obtained in $O(n^5)$ time and it follows that the algorithm 
can be implemented to run in $O(n^6)$ time. 

\smallskip 

There exists a black-and-white coloring of $G$ 
with $b$ black vertices and $w$ white vertices if and only 
if $\gamma_1(b,w) = \text{true}$. This proves the theorem.     
\qed\end{proof}

\section{Strongly chordal graphs}
\label{section SC}

The class of interval graphs is contained in the class of 
strongly chordal graphs. In this section we generalize the 
results of Section~\ref{section interval} to the class of 
strongly chordal graphs. 

\begin{definition}
Let $C=[x_1,\ldots,x_{2k}]$ be a cycle of even length. A 
chord $(x_i,x_j)$ in $C$ is an odd chord if the distance 
in $C$ between $x_i$ and $x_j$ is odd. 
\end{definition}

Recall that a graph is chordal if it has no induced 
cycle of length more than three~\cite{kn:dirac,kn:hajnal}.

\begin{definition}[\cite{kn:farber}]
A graph $G$ is strongly chordal if $G$ is chordal 
and each cycle in $G$ of 
even length at least six has an odd chord.
\end{definition}
 
Farber discovered the strongly chordal graphs as a 
subclass of chordal graph for which the weighted domination 
problem is polynomial. The class of graphs is closely related 
to the class of chordal bipartite graphs~\cite{kn:brouwer2}.  

\medskip 

There are many ways to characterize strongly chordal graphs.   
For example, a graph is strongly chordal if and only if 
its closed neighborhood matrix, or also, its clique matrix, is totally 
balanced%
~\cite{kn:anstee,kn:anstee2,kn:brouwer2,kn:farber,kn:hoffman,kn:lehel}.  
Strongly chordal graphs are also characterized 
by the property that they have no induced cycles of 
length more than three and 
no induced suns~\cite{kn:brouwer2,kn:farber}. 
For $k \geq 3$, a $k$-sun consists of a clique 
$C=\{c_1,\ldots,c_k\}$ and an independent set $S=\{s_1,\ldots,s_k\}$. 
Each vertex $s_i$, with $1 \leq i < k$, 
is adjacent to $c_i$ and to $c_{i+1}$ 
and $s_k$ is adjacent to $c_k$ and $c_1$. 

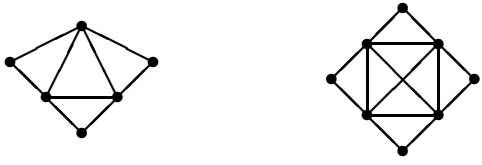
\begin{figure}
\setlength{\unitlength}{1.35pt}
\begin{center}
\begin{picture}(130,40)
\thicklines
\put(0,25){\circle*{3.0}}
\put(10,15){\circle*{3.0}}
\put(20,5){\circle*{3.0}}
\put(20,35){\circle*{3.0}}
\put(30,15){\circle*{3.0}}
\put(40,25){\circle*{3.0}}
\put(0,25){\line(2,1){20}}
\put(0,25){\line(1,-1){10}}
\put(10,15){\line(1,2){10}}
\put(10,15){\line(1,0){20}}
\put(10,15){\line(1,-1){10}}
\put(20,5){\line(1,1){10}}
\put(30,15){\line(-1,2){10}}
\put(30,15){\line(1,1){10}}
\put(40,25){\line(-2,1){20}}

\put(90,20){\circle*{3.0}}
\put(100,10){\circle*{3.0}}
\put(100,30){\circle*{3.0}}
\put(110,0){\circle*{3.0}}
\put(110,40){\circle*{3.0}}
\put(120,10){\circle*{3.0}}
\put(120,30){\circle*{3.0}}
\put(130,20){\circle*{3.0}}
\put(90,20){\line(1,-1){10}}
\put(90,20){\line(1,1){10}}
\put(100,10){\line(1,-1){10}}
\put(100,10){\line(1,0){20}}
\put(100,10){\line(1,1){20}}
\put(100,10){\line(0,1){20}}
\put(100,30){\line(1,0){20}}
\put(100,30){\line(1,1){10}}
\put(100,30){\line(1,-1){20}}
\put(110,0){\line(1,1){10}}
\put(120,10){\line(0,1){20}}
\put(120,10){\line(1,1){10}}
\put(120,30){\line(1,-1){10}}
\put(120,30){\line(-1,1){10}}
\end{picture}
\caption{A chordal graph is strongly chordal if it has no sun. 
The figure shows a 3-sun and a 4-sun.} 
\label{suns}
\end{center}
\end{figure}

\medskip 

Another way to characterize strongly chordal graphs is by the 
property that every induced subgraph has a simple vertex. 

\begin{definition}
A vertex $x$ in a graph $G$ is simple if for all 
$y,z \in N[x]$ 
\[N[y] \subseteq N[z] \quad\text{or}\quad N[z] \subseteq N[y].\]
\end{definition}
Notice that a simple vertex is simplicial, that is, its neighborhood 
is a clique. 

\begin{theorem}[\cite{kn:brouwer,kn:farber}]
A graph is strongly chordal if and only if every 
induced subgraph has a simple vertex. 
\end{theorem}

\subsection{Strongly chordal $\mathbf{k}$-trees}

We use Lehel's decomposition which decomposes a strongly 
chordal graph $G$ into a sequence of strongly chordal $k$-trees, 
for $k=1,\ldots,n-1$ such that every maximal clique of $G$ 
with cardinality $\ell+1$ is a maximal clique in the 
strongly chordal $\ell$-tree. 

\medskip 

In this section we show that there is a polynomial-time algorithm 
which solves the black-and-white coloring problem on 
strongly chordal $k$-trees. 

\begin{definition}[\cite{kn:beineke,kn:kloks,kn:moon,kn:rose}]
A $k$-tree is a connected chordal graph which is either 
a $k$-clique or in which in which every 
maximal clique has cardinality $k+1$. 
\end{definition}

Let $H=(V,E)$ be a strongly chordal $k$-tree and let $x \in V$. 
Let $C$ be a component of $H-N[x]$ and let 
\[S=N(C)=\{x_1,\ldots,x_k\}.\] 
Since $H$ has no sun, the neighborhoods in $C$ of any two 
vertices $x_i$ and $x_j$ are comparable. We assume 
that the vertices of $S$ are ordered such that 
\[N(x_1) \cap C \;\subseteq\; \ldots \;\subseteq\; N(x_k) \cap C.\]

\medskip 

The component $C$ has exactly one vertex $c$ which is adjacent to 
all vertices of $S$. Let $C_1,\ldots,C_t$ be the components 
of $H[C]-c$. For each component $C_i$ we have that 
\[N(C_i)=\{c\} \cup \{x_2,\ldots,x_k\}.\] 
Obviously, the neighborhoods in $C_i$ of $x_2,\ldots,x_k$ 
are ordered by inclusion as above and there exists 
some function $f:\{1,\ldots,t\} \rightarrow \{1,\ldots,k\}$ 
such that 
\[\forall_{2 \leq \ell \leq k}\; 
N(x_{\ell}) \cap C_i \begin{cases} 
\subseteq N(c) \cap C_i & \text{if $\ell \leq f(i)$, and} \\
\supseteq N(c) \cap C_i & \text{if $\ell > f(i)$.}
\end{cases}\]

\medskip 

Consider a black-and-white coloring of the vertices 
in $C$ with $b$ black vertices and $w$ white vertices. 
Assume that the vertex $c$ is colored black. 
For $i=1,\ldots,t$, 
let $\gamma_i(b_i,w_i)=\text{true}$ if there exists a 
black-and-white coloring of $H[C_i]$ with $b_i$ black vertices 
and $w_i$ white vertices such that $c$ is not adjacent to 
any white vertex in $C_i$. 
Let $\gamma(b,w)=\text{true}$ if there exists a 
black-and-white coloring of the vertices in $C$ 
such that $c$ is colored black. 
Then 
\begin{eqnarray*}
\gamma(b,w)=\text{true} && \quad\text{if and only if}\quad\\
&& \exists_{b_1} \cdots \exists_{b_t}\;  
\exists_{w_1} \cdots \exists_{w_t} \\
&& b=1+\sum_{i=1}^t b_i \quad\text{and}\quad 
w=\sum_{i=1}^t w_i \quad\text{and}\\
&& \forall_{1 \leq i \leq t} \; \gamma_i(b_i,w_i) = \text{true}.    
\end{eqnarray*}

\medskip 

Similar formulas can be obtained for the cases where 
$c$ is colored white and where $c$ is uncolored. 

\medskip 

Notice that, in order to maintain $\gamma_i(b_i,w_i)$, 
it is sufficient to keep a table for each possible position  
$f(i)$ which the vertex $c$ can occupy in the 
neighborhood ordering of the component $C_i$. 

\medskip 

Obviously, when $c$ is colored black, no vertex of 
$S$ can be colored white, since $c$ is 
adjacent to all vertices of $S$. Assume that a vertex $x_{\ell} \in S$ 
is colored black. Then $x_{\ell}$ is not adjacent to any  
white vertex in $C$. 
In order to know whether we can color $x_{\ell}$ black, 
it is sufficient to have the index in the neighborhood 
ordering of each $N(C_i)$, of the vertex (if any) with the 
smallest neighborhood in $C_i$ which is 
adjacent to any white vertex. 

\medskip 

Summarizing, it suffices to keep a table 
of boolean values $\gamma(b,w,p,q)$. The value of 
$\gamma(b,w,p,q)$ is true 
if there exists a black-and-white coloring 
of $H[C]$ with $b$ black vertices and $w$ white vertices 
such that 
\begin{enumerate}[\rm 1.]
\item $x_1,\ldots,x_p$ is not adjacent to any black vertex and, 
if $p < k$,  
$x_{p+1}$ is adjacent to a black vertex, and 
\item $x_1,\ldots,x_q$ is not adjacent to any white vertex and, 
if $q < k$,  
$x_{q+1}$ is adjacent to a white vertex. 
\end{enumerate}
Notice that, {\em e.g.\/}, $p=0$ implies 
that $x_1$ is adjacent to a black vertex in $C$, 
that is, the vertex 
$c$ is colored black. 

We omit the lengthy description 
of the recursive formula for $\gamma(b,w,p,q)$. 

\medskip 

Consider a vertex $x$. Possibly, $H-N[x]$ contains more than 
one component. In order to combine the black-and-white 
colorings of different components in a table, we proceed as 
follows. 

\medskip 

Choose a simplicial vertex $r$ in $H$ a `root.' Let $x$ be a 
vertex which is not adjacent to $r$. 
Let $C_x(r)$ 
be the component of $H-N[x]$ that contains $r$. 
Let the vertices in $N(C_x(r))$ be ordered  
\begin{eqnarray*}
N(C_x(r))&=&\{\;y_1,\;\ldots,\;y_k\;\} \quad\text{such that}\\
&& N(y_1) \cap C_x(r) \; \subseteq \;\ldots\;\subseteq\;N(y_k) \cap C_x(r).
\end{eqnarray*}
 
{F}rom now on, we consider only pairs $x$ and $C$ such that 
$x$ is not adjacent to $r$ and such that $C$ is a component of 
$G-N[x]$ that does not contain $r$. 
For such a pair consider the vertex $c \in C$ which is 
adjacent to all vertices of $S=N(C)$. 
Let $S=\{y_1,\ldots,y_k\}$ be the vertices of $S$ 
ordered such that 
\[N(y_1) \cap C_c(r) \;\subseteq\; \ldots\;\subseteq\; 
N(y_k) \cap C_c(r).\] 

\medskip 

Suppose we want to color a vertex $y \in C_x(r)$ black. 
Let $\ell$ be the smallest index such that $y_{\ell}$ 
is adjacent to $y$. We need to make sure that no 
vertex of $\{y_{\ell},\ldots,y_k\}$ is colored black. 

\medskip 

For that purpose, 
define the boolean variable $\gamma^{\prime}(b,w,p,q)$ 
as true if there exists a black-and-white coloring 
of the vertices in $C\cup S$ with $b$ black vertices 
and $w$ white vertices such that 
\begin{enumerate}[\rm i.]
\item $y_p, \ldots,y_k$ are not colored black, and 
\item $y_q, \ldots, y_k$ are not colored white. 
\end{enumerate}
Notice that the values $\gamma^{\prime}(b,w,p,q)$ can easily be 
deduced from the $\gamma$-table(s). 

\begin{theorem}
There exists a polynomial-time algorithm which solves 
the black-and-white coloring problem on $k$-trees.
\end{theorem}
\begin{proof}
The algorithm sorts the pairs $(x,C)$, where $x \in V$ 
is a vertex not adjacent to $r$ and where $C$ is a 
component of $H-N[x]$ that does not contain $r$, 
in increasing order of $|C|$. 
There are $O(n^2)$ such pairs since each pair $(x,C)$ is fixed 
by the pair of vertices $x$ and $c$, where $c \in C$ is 
the vertex in $C$ with 
\[N(C)=N(x) \cap N(c).\] 
For each pair 
compute a table of boolean values $\gamma(b,w,p,q)$ from the 
tables at the components $C_1,\ldots,C_t$ as outlined   
above. The components are added one by one. 
When a component $C_i$ is handled an update 
is made for the suitable table entries of $C$ by going through the 
table entries $\gamma_i(b_i,w_i,p_i,q_i)$ of the component $C_i$. 
There are $O(n^2k^2)$ entries in each table, and since 
$t \leq n$, a table for $(x,C)$ is computed in $O(n^3k^2)$ 
time. 

\smallskip 

In a similar manner compute the tables with boolean 
values $\gamma^{\prime}(b,w,p,q)$. As above, it 
is easy to update a table 
for a vertex $x$ when there are two or more 
components of $G-N[x]$ that do not contain $r$.  
\qed\end{proof}

\subsection{The transition from $\mathbf{k}$-trees to 
$\mathbf{(k+1)}$-trees}

Lehel's decomposition for strongly chordal graphs $G=(V,E)$ 
is a sequence of strongly chordal $k$-trees with vertex set $V$ for 
$k=1,\ldots,n-1$ such that every maximal clique of $G$ is a 
maximal clique in one of the strongly chordal $k$-trees. 
The $(k+1)$-tree in this sequence is obtained from the 
$k$-tree by a construction that we describe next. 

\medskip 

Let $x$ and $y$ be nonadjacent vertices in a graph $G$. 
An $x,y$-separator is a set of vertices $S$ such that 
$x$ and $y$ are in different components of $G-S$. 
The $x,y$-separator is minimal if no proper subset of 
$S$ separates $x$ and $y$ in different components. 
A set $S$ is a minimal separator in $G$ if there exist 
nonadjacent vertices $x$ and $y$ such that $S$ is a minimal 
$x,y$-separator. 

Rose characterizes chordal graphs by the property that 
every minimal separator is a clique~\cite{kn:rose2}. 
In a $k$-tree $H$ every minimal separator is a $k$-clique~\cite{kn:rose}. 
Consider a pair $(x,C)$ where $x$ is a vertex in $H$ 
and where $C$ is a component of $H-N[x]$. Then $N(C)$ is a minimal 
separator since it separates $x$ from every vertex in $C$. If 
$c \in C$ is adjacent to all vertices of $N(C)$ then $N(C)$ is 
the common neighborhood of $x$ and $c$ and so, no proper subset 
of $N(C)$ separates $x$ and $c$.  
Furthermore, it is easy to see 
that every minimal separator in a $k$-tree is of this form. 

\medskip 

Let $T_k$ be a clique tree for a $k$-tree $H$. 
A clique tree $T_k$ 
for $H$ is a tree of which the vertices are the maximal 
cliques in $H$. The tree $T_k$ satisfies the following 
property. 
\begin{quote}
For every vertex $x$ in $H$ the maximal cliques 
that contain $x$ form a subtree of $T_k$. 
\end{quote}

\medskip 

Notice that, if $C_1$ and $C_2$ are adjacent cliques in 
$T_k$ then $C_1 \cap C_2$ is a minimal separator in $H$.  
Since every minimal separator in $H$ is a $k$-clique,    
\[|C_1 \cap C_2|=k.\] 
Thus each edge in $T_k$ corresponds with a minimal 
separator in $H$ and it is easy to see that this collection of 
minimal separators is the set of all the minimal separators in $H$. 

\medskip 

Define a $(k+1)$-tree $H^{\prime}$ as follows~\cite{kn:lehel}. 
The maximal 
cliques of $H^{\prime}$ are the unions of maximal cliques 
that are endpoints of edges in $T_k$. By the observation 
above, these maximal cliques have cardinality $k+2$. 
A clique tree $T_{k+1}$ for $H^{\prime}$ is a spanning tree 
of the linegraph $L(T_k)$. 

\medskip 

Lehel proves that for every strongly chordal graph $G$ there 
is a sequence $H_k$ of $k$-trees such that every 
maximal clique in $G$ is a maximal clique in one of 
the $H_k$. Furthermore, each $H_{k+1}$ with 
clique tree $T_{k+1}$ is obtained from $H_k$ 
with clique tree $T_k$  
by an operation  
as described above~\cite{kn:lehel}. The starting 
clique tree $T_0$ is called the basic tree in~\cite{kn:lehel}. 
This basic tree $T_0$ is any tree with vertex set $V$ such that every 
maximal clique in $G$ induces a subtree of $T_0$.      

\medskip 

We shortly analyze the transition of the 
$k$-tree $H$ into the $(k+1)$-tree $H^{\prime}$. 

Let $H$ be a $k$-tree with clique tree $T_k$. 
Let $x$ be a vertex in $H$ and let $C$ be a component 
of $H-N_H[x]$. Let $S=N(C)=\{x_1,\ldots,x_k\}$. 
We assume that 
\[N_H(x_1) \cap C \;\subseteq \; \ldots\;\subseteq\; N_H(x_k) \cap C.\] 

\medskip 

Let $c$ be the unique vertex 
in $C$ which is adjacent to $S$ in $H$. 
Let $C_1,\ldots,C_t$ be the components of $H[C]-c$.  
Each component $C_i$ has a unique vertex $c_i$ 
such that 
\[N_H(C_i)=N_H(c_i) \cap N_H(x_1) = \{c\} \cup \{x_2,\ldots,x_k\}.\] 

\medskip 

Assume that $S \cup \{x\}$ is the parent of $S \cup \{c\}$ 
in $T_k$. Every edge in $T_k$ merges into one clique 
of $H^{\prime}$. Thus the two $(k+1)$-cliques 
\[S \cup \{x\} \quad\text{and}\quad S \cup \{c\} \quad\text{merge into}
\quad 
S \cup \{x,c\} \quad\text{in $H^{\prime}$.}\]

\medskip 

The component $C$ that contains $c$ consists of 
vertices that appear in maximal cliques in the subtree 
of $S\cup \{c\}$.\footnote{Possibly, when $G$ 
is disconnected, the subtree contains also the vertices 
of some other components of $G$. Notice that a clique tree 
for a chordal graph may connect the clique trees of its 
components in some arbitrary way.}
Consider all the maximal cliques in $T_k$ 
that contain $\{c,x_2,\ldots,x_k\}$. Notice that this 
includes the $(k+1)$-cliques 
\[\{c_i\} \cup \{c,x_2,\ldots,x_k\}.\]
By the Helly property (see, {\em e.g.\/},~\cite{kn:gavril,kn:kloks}), 
these maximal cliques form a subtree of $T_k$ rooted at 
$S \cup \{c\}$.  
It follows that 
Lehel's construction of the $(k+1)$-tree $H^{\prime}$ 
creates a tree $R$, with vertex set  
\[\{x_1,c_1,\ldots,c_t\},\] rooted at $x_1$.
Each edge in $R$  
forms a $k+2$-clique with $S$ in $H^{\prime}$. 

\medskip 

Obviously, not every maximal clique in $H^{\prime}$ is a 
clique in $G$. 
Assume that $x$ and $c$ are not adjacent in $G$. 
Consider a vertex $y$ in $H^{\prime}$ which is adjacent to 
$S \cup \{x\}$ and which is not in $C$. Consider 
the computation of the table 
for the pair $y$ and $C$ in $H^{\prime}$. 
The vertex $x$ is not adjacent to $c$ and so it is 
not adjacent to any vertex in $C$ in $H^{\prime}$. 
In that case the variables $p$ and $q$ in $\gamma(b,w,p,q)$ 
are at least one, since $x$ is the smallest vertex in the 
neighborhood ordering and $x$ is not adjacent to any 
black or white vertex in $C$. 
Since $x$ enters the separator as a minimal vertex in the 
neighborhood ordering, the table for the pair $y$ and $C$ 
can be determined in the same manner as described in the 
previous section. 
      
\begin{theorem}
There exists a polynomial-time algorithm which solves the 
black-and-white coloring problem on strongly 
chordal graphs. 
\end{theorem}
\begin{proof}
An analysis of Lehel's decomposition of the strongly chordal graph 
into a sequence of $k$-trees shows that this 
decomposition can be obtained $O(n^4)$ 
time~\cite{kn:lehel}. 
By the result of the previous section and the 
observations above, 
the tables for each $k$-tree can be obtained in 
$O(n^5 k^2)=O(n^7)$ 
time. Since the list contains at most $n$ $k$-trees this proves 
the theorem. 
\qed\end{proof}
 
\section{Splitgraphs}
\label{section split}

In this section we show that the black-and-white coloring 
problem on splitgraphs is NP-complete. 

\begin{definition}
A graph $G=(V,E)$ is a splitgraph if there exists a 
partition of the vertices in two sets $C$ and $S$ such that 
$C$ induces a clique in $G$ and $S$ induces an independent set 
in $G$. Here, one of the two sets $C$ and $S$ may be empty. 
\end{definition}

A splitgraph can be characterized in various ways. Notice that, 
if $G$ is a splitgraph then $G$ is chordal 
and, furthermore,  its complement $\Bar{G}$ is also a 
splitgraph. Actually, this property characterizes 
splitgraphs~\cite{kn:foldes}; a graph $G$ is a 
splitgraphs if and only if $G$ and its complement 
$\Bar{G}$ are both chordal. Splitgraphs are exactly the graphs 
that have no induced $C_4$, $C_5$ or $2K_2$~\cite{kn:foldes}.        

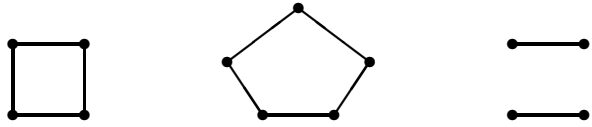
\begin{figure}
\setlength{\unitlength}{1.35pt}
\begin{center}
\begin{picture}(160,30)
\thicklines
\put(0,0){\circle*{3.0}}
\put(20,0){\circle*{3.0}}
\put(0,20){\circle*{3.0}}
\put(20,20){\circle*{3.0}}
\put(0,0){\line(1,0){20}}
\put(0,0){\line(0,1){20}}
\put(20,20){\line(-1,0){20}}
\put(20,20){\line(0,-1){20}}

\put(60,15){\circle*{3.0}}
\put(70,0){\circle*{3.0}}
\put(80,30){\circle*{3.0}}
\put(90,0){\circle*{3.0}}
\put(100,15){\circle*{3.0}}
\put(60,15){\line(2,-3){10}}
\put(60,15){\line(4,3){20}}
\put(70,0){\line(1,0){20}}
\put(90,0){\line(2,3){10}}
\put(100,15){\line(-4,3){20}}

\put(140,0){\circle*{3.0}}
\put(160,0){\circle*{3.0}}
\put(140,20){\circle*{3.0}}
\put(160,20){\circle*{3.0}}
\put(140,0){\line(1,0){20}}
\put(140,20){\line(1,0){20}}
\end{picture}
\caption{A graph is a splitgraph if it has no $C_4$, $C_5$ or $2K_2$.}
\label{split}
\end{center}
\end{figure}

\begin{theorem}
The black-and-white coloring problem is NP-complete 
for the class of splitgraphs. 
\end{theorem}
\begin{proof}
Since splitgraphs are closed under complementation, 
we can formulate the problem as a black-and-white coloring 
problem with all black vertices adjacent to all white vertices. 
We call this the `inverse B\&W-coloring problem.' 

\smallskip 

We adapt a proof of Johnson, which proves the NP-completeness 
of finding a balanced complete bipartite subgraph in a 
bipartite graph~\cite[Page~446]{kn:johnson}. 

\smallskip 

Let $G=(V,E)$ be a graph with $|V|=n$. 
Construct a splitgraph $H$ as follows. 
The clique of the splitgraph consists of the set $V$. 
The independent set of the splitgraph consists of 
the set $E$. 
In the splitgraph, make a vertex $x\in V$ adjacent to 
an edge $\{y,z\} \in E$ if and only if $x$ is {\sc not} 
an endpoint of $\{y,z\}$. 

This completes the description of $H$. 

\smallskip 

Assume that the clique number of $G$ is $\omega$. We may 
assume that $n$ is even and $n > 6$, and 
that $\omega=\frac{n}{2}$~\cite{kn:johnson}. 
 
Then we have an inverse B\&W-coloring of 
$H$ with 
\begin{equation}
\label{bw-eqn}
b=\omega \quad\text{and}\quad w=\omega+\binom{\omega}{2}=
\binom{\omega+1}{2}.
\end{equation}

\smallskip 

For the converse, assume that $H$ has an inverse  
B\&W-coloring with the numbers of black and white vertices as in 
Equation~\ref{bw-eqn}. 
Since $E$ is an independent set in $H$ the colored vertices in $E$ 
must all have the same color. First assume that $E$ contains 
no white vertices. Then $V$ contains a set 
$W$ of white vertices, and $V\setminus W$ is black. 
Since 
\[w=\omega+\binom{\omega}{2} > n=2 \omega \quad\text{if $n > 6$},\]
this is not possible. Thus the inverse black-and-white coloring 
has white vertices in $E$. 

\smallskip 

Assume that the inverse B\&W-coloring has 
a set $E^{\prime}$ of white vertices in $E$ and a set of 
$V^{\prime}$ of $\omega$ black vertices in $V$. By the construction, 
no edge of $E^{\prime}$ has an endpoint in $V^{\prime}$. Now 
$|V \setminus V^{\prime}|=\omega$ and all the endpoints 
of $E^{\prime}$ are in $V\setminus V^{\prime}$. The only 
possibility is that $E^{\prime}$ is the set of edges of 
a clique $V \setminus V^{\prime}$ of cardinality $\omega$ in $G$.  

This proves the theorem.
\qed\end{proof}


\begin{thebibliography}{99}

\bibitem{kn:acharya}Acharya,~B. and M.~Las~Vergnas, 
Hypergraphs with cyclomatic number zero, triangulated graphs, 
and an inequality, 
{\em Journal of Combinatorial Theory, Series B\/} {\bf 33} (1982), 
pp.~52--56. 

\bibitem{kn:anstee}Anstee,~R., 
Hypergraphs with no special cycles, 
{\em Combinatorica\/} {\bf 3} (1983), pp.~141--146.  

\bibitem{kn:anstee2}Anstee,~R. and M.~Farber, 
Characterizations of totally balanced matrices, 
{\em Journal of Algorithms\/} {\bf 5} (1984), pp.~215--230. 

\bibitem{kn:bandelt}Bandelt,~H. and H.~Mulder, 
Distance-hereditary graphs, 
{\em Journal of Combinatorial Theory, Series B\/} {\bf 41} (1986), 
pp.~182--208. 

\bibitem{kn:beineke}Beineke,~L. and R.~Pippert, 
The number of labeled $k$-dimensional trees, 
{\em Journal of Combinatorial Theory\/} {\bf 6} (1969), 
pp.~200-205. 

\bibitem{kn:berend}Berend,~D. and S.~Zucker, 
The black-and-white coloring problem on trees, 
{\em Journal of Graph Algorithms and Applications\/} {\bf 13} 
(2009), pp.~133--152. 

\bibitem{kn:booth}Booth,~K. and G.~Lueker, 
Linear algorithms to recognize interval graphs and test for 
the consecutive ones property, 
{\em Proceedings STOC'75\/}, ACM (1975), pp.~255--265. 

\bibitem{kn:broersma}Broersma,~H., T.~Kloks, D.~Kratsch and 
H.~M\"uller, 
Independent sets in asteroidal triple-free graphs, 
{\em SIAM Journal on Discrete Mathematics\/} {\bf 12} (1999), 
pp.~276--287. 

\bibitem{kn:brouwer2}Brouwer,~A., P.~Duchet and A.~Schrijver, 
Graphs whose neighborhoods have no special cycle, 
{\em Discrete Mathematics\/} {\bf 47} (1983), pp.~177--182. 

\bibitem{kn:brouwer}Brouwer,~A. and A.~Kolen, 
A super-balanced hypergraph has a nest point. 
Technical Report ZW~146, Mathematisch Centrum, Amsterdam, 1980. 

\bibitem{kn:chang}Chang,~M., S.~Hsieh and G.~Chen, 
Dynamic programming on distance-hereditary graphs, 
{\em Proceedings ISAAC'97\/}, Springer, LNCS~1350 (1997), pp.~344--353. 

\bibitem{kn:chvatal}Chv\'atal,~V. and P.~Hammer, 
Aggregation of inequalities in integer programming. 
Technical Report STAN-CS-75-518, Stanford University, California, 1975. 

\bibitem{kn:corneil}Corneil,~D., H.~Lerchs and L.~Stewart-Burlingham, 
Complement reducible graphs, 
{\em Discrete Applied Mathematics\/} {\bf 3} (1981), pp.~163--174. 

\bibitem{kn:corneil2}Corneil,~D., Y.~Perl and L.~Stewart, 
A linear recognition algorithm for cographs, 
{\em SIAM Journal on Computing\/} {\bf 14} (1985), pp.~926--934. 

\bibitem{kn:dirac}Dirac,~G.,
On rigid circuit graphs,
{\em Abhandlungen aus dem Mathematischen Seminar der Universit\"at
Hamburg\/} {\bf 25} (1961), pp.~71--76.

\bibitem{kn:farber}Farber,~M., 
Characterizations of strongly chordal graphs, 
{\em Discrete Mathematics\/} {\bf 43} (1983), pp.~173--189. 

\bibitem{kn:foldes}F\"oldes,~S. and P.~Hammer, 
Split graphs, 
{\em Congressus Numerantium\/} {\bf 19} (1977), pp.~311--315. 

\bibitem{kn:gavril}Gavril,~F., 
The intersection graphs of subtrees in trees are exactly the 
chordal graphs, 
{\em Journal of Combinatorial Theory, Series B\/} {\bf 16} 
(1974), pp.~47--56. 

\bibitem{kn:gilmore}Gilmore,~P. and A.~Hoffman, 
A characterization of comparability graphs and of 
interval graphs, 
{\em The Canadian Journal of Mathematics\/} {\bf 16} (1964), 
pp.~539--548. 

\bibitem{kn:hajnal}Hajnal,~A. and J.~Sur\'anyi, 
\"Uber die Aufl\"osung von Graphen in vollst\"andige Teilgraphen, 
{\em Annales Universitatis Scientiarum Budapestinensis de Rolando 
E\"otv\"os Nominatae -- Sectio Mathematicae\/} {\bf 1} (1958), 
pp.~113--121.  

\bibitem{kn:hammer}Hammer,~P., U.~Peled and X.~Sun, 
Difference graphs, 
{\em Discrete Applied Mathematics\/} {\bf 28} (1990), pp.~35--44. 

\bibitem{kn:hansen}Hansen,~P., A.~Hertz and N.~Quinodos, 
Splitting trees, 
{\em Discrete Mathematics\/} {\bf 165} (1997), pp.~403--419. 

\bibitem{kn:hoffman}Hoffman,~A., A.~Kolen and M.~Sakarovitch, 
Totally-balanced and greedy matrices. Technical Report BW~165/82, 
Mathematisch Centrum, Amsterdam, 1982. 
 
\bibitem{kn:howorka}Howorka,~E., 
A characterization of distance-hereditary graphs, 
{\em The Quarterly Journal of Mathematics\/} {\bf 28} (1977), pp.~417--420. 

\bibitem{kn:johnson}Johnson,~D., 
The NP-completeness column: An ongoing guide, 
{\em Journal of Algorithms\/} {\bf 8} (1987), pp.~438--448. 

\bibitem{kn:kloks}Kloks,~T.,
{\em Treewidth -- Computations and Approximations\/},
Springer, Lecture Notes in Computer Science {\bf 842}, 1994.

\bibitem{kn:lehel}Lehel,~J., 
A characterization of totally balanced matrices, 
{\em Discrete Mathematics\/} {\bf 57} (1985), pp.~59--65.

\bibitem{kn:lekkerkerker}Lekkerkerker,~C. and D.~Boland, 
Representation of finite graphs by a set of intervals on the 
real line, 
{\em Fundamenta Mathematicae\/} {\bf 51} (1962), pp.~45--64. 

\bibitem{kn:mahadev}Mahadev,~N. and U.~Peled, 
{\em Threshold graphs and related topics\/}, 
Elsevier Series {\em Annals of Discrete Mathematics\/} {\bf 56}, 
1995. 

\bibitem{kn:moon}Moon,~J., 
The number of labeled $k$-trees, 
{\em Journal of Combinatorial Theory\/} {\bf 6} (1969), 
pp.~196--199. 

\bibitem{kn:rose2}Rose,~D., 
Triangulated graphs and the elimination process, 
{\em Journal of Mathematical Analysis and Applications\/} {\bf 32} 
(1970), pp.~597--609. 

\bibitem{kn:rose}Rose,~D., 
On simple characterizations of $k$-trees, 
{\em Discrete Mathematics\/} {\bf 7} (1974), pp.~317--322. 

\bibitem{kn:yannakakis}Yannakakis,~M., 
The complexity of the partial order dimension problem, 
{\em SIAM Journal on Algebraic and Discrete Methods\/} {\bf 3} 
(1982), pp.~351--358. 

\end{thebibliography}
\end{document}